\input amstex
\documentstyle{amsppt}
%
\catcode`@=11
\redefine\output@{%
  \def\break{\penalty-\@M}\let\par\endgraf
  \ifodd\pageno\global\hoffset=105pt\else\global\hoffset=8pt\fi  
  \shipout\vbox{%
    \ifplain@
      \let\makeheadline\relax \let\makefootline\relax
    \else
      \iffirstpage@ \global\firstpage@false
        \let\rightheadline\frheadline
        \let\leftheadline\flheadline
      \else
        \ifrunheads@ 
        \else \let\makeheadline\relax
        \fi
      \fi
    \fi
    \makeheadline \pagebody \makefootline}%
  \advancepageno \ifnum\outputpenalty>-\@MM\else\dosupereject\fi
}
\def\Beta{\mathchar"0\hexnumber@\rmfam 42}
\catcode`\@=\active
\nopagenumbers
\chardef\textvolna='176

\chardef\bigalpha='013
\def\negskp{\hskip -2pt}

\def\compos{\,\raise 1pt\hbox{$\sssize\circ$} \,}

\chardef\degree="5E
\def\tr{\operatorname{tr}}
\accentedsymbol\tildepsi{\kern 3.2pt\tilde{\kern -3.2pt\psi}}
\accentedsymbol\ssizetildepsi{\ssize\kern 2.5pt\tilde{\kern -2.5pt\psi}}
\def\Pover#1{\overset{\kern 1.5pt #1}\to P}
\def\Qover#1{\overset{\kern 1.5pt #1}\to Q}
\def\Rover#1{\overset{\kern 1.5pt #1}\to R}
\def\Sover#1{\overset{\kern 1.5pt #1}\to S}
\def\blue#1{#1}

\catcode`#=11\def\diez{#}\catcode`#=6
\catcode`&=11\catcode`&=4
\catcode`_=11\def\podcherkivanie{_}\catcode`_=8
\catcode`~=11\catcode`~=\active
\def\mycite#1{\cite{\blue{#1}}\immediate\special{ps:
     ShrHPSdict begin /ShrBORDERthickness 0 def}}
\def\myciterange#1#2#3#4{\cite{\blue{#2#3#4}}\immediate\special{ps:
     ShrHPSdict begin /ShrBORDERthickness 0 def}}
\def\mytag#1{%
    \tag#1}
\def\mythetag#1{\thetag{\blue{#1}}\immediate\special{ps:
     ShrHPSdict begin /ShrBORDERthickness 0 def}}
\def\myrefno#1{\no#1}
\def\myhref#1#2{\blue{#2}\immediate\special{ps:
     ShrHPSdict begin /ShrBORDERthickness 0 def}}
\def\myEarXivlink{\myhref{http://arXiv.org}{http:/\negskp/arXiv.org}}

\def\mytheorem#1{\csname proclaim\endcsname{Theorem #1}}
\def\mytheoremwithtitle#1#2{\csname proclaim\endcsname{Theorem #1#2}}
\def\mythetheorem#1{\blue{#1}\immediate\special{ps:
     ShrHPSdict begin /ShrBORDERthickness 0 def}}
\def\mylemma#1{\csname proclaim\endcsname{Lemma #1}}
\def\mylemmawithtitle#1#2{\csname proclaim\endcsname{Lemma #1#2}}
\def\mythelemma#1{\blue{#1}\immediate\special{ps:
     ShrHPSdict begin /ShrBORDERthickness 0 def}}
\def\mycorollary#1{\csname proclaim\endcsname{Corollary #1}}

\def\mydefinition#1{\definition{Definition #1}}
\def\mythedefinition#1{\blue{#1}\immediate\special{ps:
     ShrHPSdict begin /ShrBORDERthickness 0 def}}
\def\myconjecture#1{\csname proclaim\endcsname{Conjecture #1}}
\def\myconjecturewithtitle#1#2{\csname proclaim\endcsname{Conjecture #1#2}}

\def\myproblem#1{\csname proclaim\endcsname{Problem #1}}
\def\myproblemwithtitle#1#2{\csname proclaim\endcsname{Problem #1#2}}


\pagewidth{360pt}
\pageheight{606pt}
\topmatter
\title
Umbilical and zero curvature equations in a class of second 
order ODE's 
\endtitle
\rightheadtext{Umbilical and zero curvature equations \dots}
\author
Ruslan Sharipov
\endauthor
\address Bashkir State University, 32 Zaki Validi street, 450074 Ufa, Russia
\endaddress
\email\myhref{mailto:r-sharipov\@mail.ru}{r-sharipov\@mail.ru}
\endemail
\abstract
     The class of second order ODE's cubic with respect to the first order 
derivative is considered. Using geometric structures associated with these 
equations, the subclasses of umbilical equations, zero mean curvature 
equations, and zero Gaussian curvature equations are defined. Zero mean 
curvature equations are studied within the framework of the first 
case of intermediate degeneration with the stress on their pseudoscalar 
and scalar invariants.
\endabstract
\endtopmatter
\loadbold
\TagsOnRight
\document

\head
1. Introduction.
\endhead
     Since the epoch of classical papers (see \mycite{1} and \mycite{2}) the class 
of second order differential equations cubic with respect to the first order derivative 
$$
\hskip -2em
y''=P(x,y)+3\,Q(x,y)\,y'+3\,R(x,y)\,(y')^2+S(x,y)\,(y')^3
\mytag{1.1}
$$
attracted the attention due to rich geometric structures associated with equations of
this class. This class is closed with respect to transformations of the form 
$$
\hskip -2em
\cases
\tilde x=\tilde x(x,y),\\
\tilde y=\tilde y(x,y),
\endcases
\mytag{1.2}
$$
which can be interpreted as changes of local curvilinear coordinates in $\Bbb R^2$ or
in some two-dimensional manifolds. About 19 years ago in \mycite{3} and \mycite{4} 
the equations \mythetag{1.1} were classified using their scalar invariants 
derived from geometric structures associated with them. They were subdivided into 
{\bf nine subclasses} closed with respect to transformations of the form 
\mythetag{1.2}. Here is the list of these classes, which are called cases
in \mycite{3} and \mycite{4}:
\roster
\item"---" the case of general position (the richest class);
\item"---" the first case of intermediate degeneration;
\item"---" the second case of intermediate degeneration;
\item"---" the third case of intermediate degeneration;
\item"---" the fourth case of intermediate degeneration;
\item"---" the fifth case of intermediate degeneration;
\item"---" the sixth case of intermediate degeneration;
\item"---" the seventh case of intermediate degeneration;
\item"---" the case of maximal degeneration (the smallest class);
\endroster
The case of general position was previously studied in \mycite{5} using the
same geometric methods as in \mycite{3} and \mycite{4}. Some time later than
\mycite{3} and \mycite{4} other approaches to classifications of the equations 
\mythetag{1.1} were considered (see \myciterange{3}{6}{--}{10}) and 
Yu\.\,Yu\.\,Bagderina's classification in \mycite{11}).\par 
     In this paper we define three geometric subclasses of the equations 
\mythetag{1.1} that can potentially intersect each of the last eight classes 
in the above classification and carefully study one of them in the first case 
of intermediate degeneration.\par
\head
2. Some notations and definitions.
\endhead
    Transformations of the form \mythetag{1.2} are interpreted as changes of
local coordinates. They are assumed to be locally invertible. The inverse 
transformations for them are written similarly in the following form:
$$
\hskip -2em
\cases
x=\tilde x(\tilde x,\tilde y),\\
y=\tilde y(\tilde x,\tilde y).
\endcases
\mytag{2.1}
$$
Like in \myciterange{3}{3}{--}{5} and \mycite{12}, here we use dot index 
notations for partial derivatives, e\.\,g\. having two functions $f(x,y)$ and  
$g(\tilde x,\tilde y)$ we write
$$
\xalignat 2
\hskip -2em
f_{\ssize p.q}=\frac{\partial^{p+q}f}{\partial x^p\,\partial y^q},
&&g_{\ssize p.q}=\frac{\partial^{p+q}g}{\partial\tilde x^p\,
\partial\tilde y^q}.
\mytag{2.2}
\endxalignat
$$
Then we write the Jacoby matrices of the direct and inverse transformations 
\mythetag{1.2} and \mythetag{2.1} in terms of the above notations 
\mythetag{2.2}:
$$
\xalignat 2
&\hskip -2em
S=\Vmatrix
x_{\sssize 1.0} &x_{\sssize 0.1}\\
\vspace{1ex}
y_{\sssize 1.0} &y_{\sssize 0.1}
\endVmatrix,
&&T=\Vmatrix
\tilde x_{\sssize 1.0} &\tilde x_{\sssize 0.1}\\
\vspace{1ex}
\tilde y_{\sssize 1.0} &\tilde y_{\sssize 0.1}
\endVmatrix.
\mytag{2.3}
\endxalignat
$$
In differential geometry the Jacoby matrices \mythetag{2.3} are called the 
direct and inverse transition matrices respectively (see \mycite{13}).\par
     Tensorial and pseudotensorial fields in local coordinates are presented
as arrays of functions whose arguments are local coordinates $x,\,y$ or 
$\tilde x,\,\tilde y$ respectively. These arrays of functions are called their 
components. They obey some definite transformation rules under a change of
local coordinates.
\mydefinition{2.1} A pseudotensorial field of the type $(r,s)$ and weight $m$ 
is an array of functions $F^{i_1\ldots\,i_r}_{j_1\ldots\,j_s}$ which under 
the change of coordinates \mythetag{1.2} transforms as 
$$
\hskip -2em
F^{i_1\ldots\,i_r}_{j_1\ldots\,j_s}=
(\det T)^m\sum\Sb p_1\ldots p_r\\ q_1\ldots q_s\endSb
S^{i_1}_{p_1}\ldots\,S^{i_r}_{p_r}\,\,
T^{q_1}_{j_1}\ldots\,T^{q_s}_{j_s}\,\,
\tilde F^{p_1\ldots\,p_r}_{q_1\ldots\,q_s}.
\mytag{2.4}
$$
Tensorial fields are those pseudotensorial fields whose weight $m$ in \mythetag{2.4}
is zero. The prefix ``pseudo'' always indicates the nonzero weight $m\neq 0$. 
\enddefinition
     Tensorial and pseudotensorial fields of the type $(1,0)$ are called vectorial 
and pseudovectorial fields. Tensorial and pseudotensorial fields of the type $(0,1)$ 
are called covectorial and pseudocovectorial fields. And finally, scalar and 
pseudoscalar fields are those fields \pagebreak whose type is $(0,0)$.
\mydefinition{2.2} Tensorial and pseudotensorial fields whose components are 
expressed through $y'$, through the coefficients $P$, $Q$, $R$, $S$ of the equation
\mythetag{1.1}, and through their partial derivatives are called tensorial and 
pseudotensorial invariants of this equation respectively.
\enddefinition
\head
3. Some basic structures.
\endhead	
     Each equation of the form \mythetag{1.1} is associated with two 
pseudocovectorial fields $\boldsymbol\alpha$ and $\boldsymbol\beta$ and
with one pseudoscalar field $F$. The components of $\boldsymbol\alpha$
are 
$$
\hskip -2em
\aligned
&\aligned
\alpha_1=A=P_{\sssize 0.2}&-2\,Q_{\sssize 1.1}+R_{\sssize 2.0}+
 2\,P\,S_{\sssize 1.0}+S\,P_{\sssize 1.0}-\\
 \vspace{0.5ex}
 &-3\,P\,R_{\sssize 0.1}-3\,R\,P_{\sssize 0.1}-
 3\,Q\,R_{\sssize 1.0}+6\,Q\,Q_{\sssize 0.1},
 \endaligned\\
 \vspace{1ex}
&\aligned
\alpha_2=B=S_{\sssize 2.0}&-2\,R_{\sssize 1.1}+Q_{\sssize 0.2}-
 2\,S\,P_{\sssize 0.1}-P\,S_{\sssize 0.1}+\\
 \vspace{0.5ex}
 &+3\,S\,Q_{\sssize 1.0}+3\,Q\,S_{\sssize 1.0}+
 3\,R\,Q_{\sssize 0.1}-6\,R\,R_{\sssize 1.0}.
 \endaligned
\endaligned
\mytag{3.1}
$$
The weight of the field $\boldsymbol\alpha$ with the components
\mythetag{3.1} is equal to $1$ (see \myciterange{3}{3}{--}{5}).  
The components of $\boldsymbol\beta$ are expressed through $A$
and $B$ taken from \mythetag{3.1}:
$$
\xalignat 2
&\hskip -2em
\beta_1=-H,
&&\beta_2=G,
\mytag{3.2}
\endxalignat
$$
where $G$ and $H$ are given by the formulas
$$
\aligned
G&=-B\,B_{\sssize 1.0}-3\,A\,B_{\sssize 0.1}+4\,B\,
A_{\sssize 0.1}+3\,S\,A^2-6\,R\,B\,A+3\,Q\,B^2,\\
\vspace{1ex}
H&=-A\,A_{\sssize 0.1}-3\,B\,A_{\sssize 1.0}+4\,A\,
B_{\sssize 1.0}-3\,P\,B^2+6\,Q\,A\,B-3\,R\,A^2.
\endaligned
\quad
\mytag{3.3}
$$
The weight of the field $\boldsymbol\beta$ with the components
\mythetag{3.2} is equal to $3$ (see \myciterange{3}{3}{--}{5}).
\par
     In $\Bbb R^2$ and in any two-dimensional manifold there are two 
pseudotensorial fields with constant components. They are denoted by 
the same symbol $\bold d$ and are given by the same skew-symmetric 
matrix in any local coordinates: 
$$
\xalignat 2
&\hskip -2em
d_{ij}=\Vmatrix\format \r&\quad\l\\ 0 & 1\\-1 & 0\endVmatrix,
&&d^{\kern 1pt ij}=\Vmatrix\format \r&\quad\l\\ 0 & 1\\-1 & 0\endVmatrix.
\mytag{3.4}
\endxalignat
$$
The weight of the first field in \mythetag{3.4} is $-1$, the weight of 
the second field is $1$. The skew-symmetric fields \mythetag{3.4} here
play the same role as metric tensors in metric geometry. They are used
for raising and lowering indices of other fields. Raising the indices
of $\boldsymbol\alpha$ and $\boldsymbol\beta$, we get two pseudovectorial
fields
$$
\xalignat 2
&\hskip -2em
\alpha^i=\sum^2_{k=1}d^{\kern 1pt ik}\,\alpha_k,
&&\beta^i=\sum^2_{k=1}d^{\kern 1pt ik}\,\beta_k.
\mytag{3.5}
\endxalignat
$$ 
The weights of the fields $\boldsymbol\alpha$ and $\boldsymbol\beta$ with
the components \mythetag{3.5} are $2$ and $4$ respectively. These 
pseudovectorial fields are denoted with the same symbols $\boldsymbol\alpha$ 
and $\boldsymbol\beta$ as the pseudocovectorial fields in \mythetag{3.1} and 
\mythetag{3.2}. The formulas \mythetag{3.5} yield
$$
\xalignat 2
&\hskip -2em
\alpha^1=B,
&&\alpha^2=-A,
\mytag{3.6}\\
\vspace{1ex}
&\hskip -2em
\beta^1=G,
&&\beta^2=H.
\mytag{3.7}
\endxalignat
$$
The pseudoscalar field $F$ mentioned above is the third field associated
with any equation of the form \mythetag{1.1}. It is expressed through
the quantities $A$, $B$, $G$, $H$ in \mythetag{3.1}, \mythetag{3.2},
\mythetag{3.3}, \mythetag{3.6}, \mythetag{3.7} by means of the formulas 
$$
\hskip -2em
3\,F^5=\sum^2_{i=1}\alpha_i\,\beta^i=-\sum^2_{i=1}\beta_i\,\alpha^i=
A\,G+B\,H.
\mytag{3.8}
$$
The weight of the field $F$ introduced by \mythetag{3.8} is equal to $1$
(see \myciterange{3}{3}{--}{5}). As for the formula \mythetag{3.8} itself,
it can be written in a more explicit form:
$$
\hskip -2em
\aligned
F=\Bigl(A\,B\,A_{\sssize 0.1}&+B\,A\,B_{\sssize 1.0}-
A^2\,B_{\sssize 0.1}-B^2\,A_{\sssize 1.0}-\\
&-P\,B^3+3\,Q\,A\,B^2-3\,R\,A^2\,B+S\,A^3\Bigr)^{1/5}.
\endaligned
\mytag{3.9}
$$
\par
     The case of {\bf general position} is introduced by the condition
$F\neq 0$ in terms of the field \mythetag{3.9} (see \myciterange{3}{3}{--}{5}). 
This condition is equivalent to  
$$
\xalignat 3
&\hskip -2em
\boldsymbol\alpha\neq 0,
&&\boldsymbol\beta\neq 0,
&&\boldsymbol\alpha\nparallel\boldsymbol\beta.
\mytag{3.10}
\endxalignat
$$
The case of {\bf maximal degeneration} is opposite to \mythetag{3.10}.
It is given by the condition $\boldsymbol\alpha=0$, which implies
$\boldsymbol\beta=0$ and $F=0$. As for the cases of {\bf intermediate
degeneration}, they are introduced by the conditions
$$
\xalignat 2
&\hskip -2em
\boldsymbol\alpha\neq 0,
&&\boldsymbol\alpha\parallel\boldsymbol\beta,
\mytag{3.11}
\endxalignat
$$
which imply $F=0$. Each particular case of of intermediate degeneration
is specified by some auxiliary conditions added to \mythetag{3.11}
(see \myciterange{3}{3}{, }{4}).\par
\head
4. Auxiliary structures common for all cases 
of intermediate degeneration. 
\endhead
    As soon as the conditions \mythetag{3.11} are fulfilled, some new 
geometric structures associated with the equations \mythetag{1.1} arise. 
They consists of two pseudoscalar field $M$ and $N$, a pseudocovectorial 
field $\boldsymbol\gamma$, and connection components $\varGamma^k_{ij}$. 
The pseudoscalar field $N$ is introduced as a factor relating two parallel 
pseudocovectorial fields one of which is nonzero. It is given by the 
following formula:
$$
\hskip -2em
\boldsymbol\beta=3\,N\,\boldsymbol\alpha.
\mytag{4.1}
$$
In terms of the components of $\boldsymbol\alpha$ and $\boldsymbol\beta$
the formula \mythetag{4.1} yields
$$
\xalignat 2
&\hskip -2em
N=\frac{G}{3\,B},&&N=-\frac{H}{3\,A}.
\mytag{4.2}
\endxalignat
$$
The first formula applies in the case $B\neq 0$, the second one in the case 
$A\neq 0$. If both $A$ and $B$ are nonzero, both formulas are applicable.
The quantities $A$ and $B$ cannot vanish simultaneously since they are
components of the field $\boldsymbol\alpha$ and $\boldsymbol\alpha\neq 0$
in all cases of intermediate degeneration (see \mythetag{3.11}). The weight 
of the field $N$ in \mythetag{4.1} and \mythetag{4.2} is equal to $2$.\par
     The pseudoscalar field $M$ is also introduced by means of two formulas
one of which is for the case $B\neq 0$ and the second one is for $A\neq 0$
(see \mycite{4}):
$$
\gather
\hskip -2em
\aligned
M=-\frac{12\,A\,N\,(A\,S-B_{\sssize 0.1})}{5\,B}
&-A\,N_{\sssize 0.1}+\frac{24}{5}\,A\,N\,R-\\
\vspace{1ex}
-\frac{6}{5}\,N\,A_{\sssize 0.1}-\frac{6}{5}\,N\,
&B_{\sssize 1.0}+B\,N_{\sssize 1.0}-\frac{12}{5}\,B\,N\,Q,
\endaligned
\mytag{4.3}\\
\vspace{2ex}
\hskip -2em
\aligned
M=-\frac{12\,B\,N\,(B\,P+A_{\sssize 1.0})}{5\,A}
&+B\,N_{\sssize 1.0}+\frac{24}{5}\,B\,N\,Q+\\
\vspace{1ex}
+\frac{6}{5}\,N\,B_{\sssize 1.0}+\frac{6}{5}\,N\,
&A_{\sssize 0.1}-A\,N_{\sssize 0.1}-\frac{12}{5}\,A\,N\,R.
\endaligned
\mytag{4.4}
\endgather
$$
The weight of the field $M$ in \mythetag{4.3} and \mythetag{4.4} 
is equal to $4$.\par
     The pseudocovectorial field $\boldsymbol\gamma$ is introduced by
two pairs of formulas for its components, one pair is for $B\neq 0$ 
and the other is for $A\neq 0$ (see \mycite{4}):
$$
\align
&\hskip -2em
\aligned
\gamma_1=&\frac{6\,A\,N\,(A\,S-B_{\sssize 0.1})}{5\,B^2}
-\frac{18\,N\,A\,R}{5\,B}+\\
\vspace{1ex}
&+\frac{6\,N\,(A_{\sssize 0.1}+B_{\sssize 1.0})}{5\,B}
-N_{\sssize 1.0}+\frac{12}{5}\,N\,Q-2\,\Omega\,A,
\endaligned
\quad
\mytag{4.5}\\
\vspace{2ex}
&\hskip -2em
\aligned
\gamma_2=-\frac{6\,N\,(A\,S-B_{\sssize 0.1})}{5\,B}
-N_{\sssize 0.1}+\frac{6}{5}\,N\,R-2\,\Omega\,B.
\endaligned
\quad
\mytag{4.6}
\endalign
$$
The formulas \mythetag{4.5} and \mythetag{4.6} are used if $B\neq 0$. 
If $A\neq 0$, we write:
$$
\align
&\hskip -2em
\aligned
\gamma_1=\frac{6\,N\,(B\,P+A_{\sssize 1.0})}{5\,A}
-N_{\sssize 1.0}-\frac{6}{5}\,N\,Q-2\,\Omega\,A,
\endaligned
\quad
\mytag{4.7}\\
\vspace{2ex}
&\hskip -2em
\aligned
\gamma_2=-&\frac{6\,B\,N\,(B\,P+A_{\sssize 1.0})}{5\,A^2}
+\frac{18\,N\,B\,Q}{5\,A}+\\
\vspace{1ex}
&+\frac{6\,N\,(B_{\sssize 1.0}+A_{\sssize 0.1})}{5\,A}
-N_{\sssize 0.1}-\frac{12}{5}\,N\,R-2\,\Omega\,B.
\endaligned
\quad
\mytag{4.8}
\endalign
$$
The weight of the field $\boldsymbol\gamma$ with the components 
\mythetag{4.5}, \mythetag{4.6}, \mythetag{4.7}, \mythetag{4.8} is $2$. 
Raising  indices in these formulas, we get a pseudovectorial field 
denoted by the same symbol:
$$
\hskip -2em
\gamma^i=\sum^2_{k=1}d^{\kern 1pt ik}\,\gamma_k.
\mytag{4.9}
$$ 
The formula \mythetag{4.9} can be written in a more explicit form:
$$
\xalignat 2
&\hskip -2em
\gamma^1=C=\gamma_2,
&&\gamma^2=D=-\gamma_1.
\mytag{4.10}
\endxalignat
$$
Here $C$ and $D$ are notations for the components of the pseudovectorial
field $\boldsymbol\gamma$. Its weight is $3$. The formula \mythetag{4.9}
is an analog of the formulas \mythetag{3.5}, while the formulas
\mythetag{4.10} are analogous to the formulas \mythetag{3.6} and 
\mythetag{3.7}.\par
     The field $M$ given by the formulas \mythetag{4.3} and \mythetag{4.4}
is related to the fields $\boldsymbol\alpha$ and $\boldsymbol\gamma$ by
means of the formulas similar to \mythetag{3.8}: 
$$
\hskip -2em
M=\sum^2_{i=1}\alpha_i\,\gamma^i=-\sum^2_{i=1}\gamma_i\,\alpha^i
=A\,C+B\,D.
\mytag{4.11}
$$
\par
     The connection components $\varGamma^k_{ij}$ constitute the fourth
auxiliary structure common for all cases of {\bf intermediate degeneration}.
They are given by the formula
$$
\hskip -2em
\varGamma^k_{ij}=\theta^k_{ij}-\frac{\varphi_i\,\delta^k_j+
\varphi_j\,\delta^k_i}{3},
\mytag{4.12}
$$
where $\delta^k_j$ and $\delta^k_i$ are Kronecker deltas. The quantities 
$\theta^k_{ij}$ in \mythetag{4.12} are given by the components of a fully
symmetric array $\theta_{ijk}$ upon raising one of its indices:
$$
\hskip -2em
\theta^k_{ij}=\sum^2_{r=1}d^{kr}\,\theta_{rij}.
\mytag{4.13}
$$
The components $\theta_{rij}$ of the array $\theta$ in \mythetag{4.13} are
given explicitly:
$$
\xalignat 2
&\hskip -2em
\theta_{111}=P,
&&\theta_{112}=\theta_{121}=\theta_{211}=Q,\\
\vspace{-1.7ex}
&&&\mytag{4.14}\\
\vspace{-1.7ex}
&\hskip -2em
\theta_{122}=\theta_{212}=\theta_{221}=R,
&&\theta_{222}=S.
\hskip -2em
\endxalignat
$$
In addition to $\theta^k_{ij}$ defined through \mythetag{4.13} and
\mythetag{4.14}, the quantities $\varphi_i$ and $\varphi_j$ are used
in \mythetag{4.12}. They are given by the formulas 
$$
\hskip -2em
\aligned
&\varphi_1=-3\,A\,\frac{A\,S-B_{\sssize 0.1}}{5\,B^2}
-3\,\frac{A_{\sssize 0.1}+B_{\sssize 1.0}-3\,A\,R}{5\,B}
-\frac{6}{5}\,Q,\\
\vspace{2ex}
&\varphi_2=3\,\frac{A\,S-B_{\sssize 0.1}}{5\,B}
-\frac{3}{5}\,R.
\endaligned
\mytag{4.15}
$$
The formulas \mythetag{4.15} apply in the case $B\neq 0$. If $A\neq 0$, we
use the formulas 
$$
\hskip -2em
\aligned
&\varphi_1=
-3\,\frac{B\,P+A_{\sssize 1.0}}{5\,A}
+\frac{3}{5}\,Q,\\
\vspace{2ex}
&\varphi_2=
3\,B\,\frac{B\,P+A_{\sssize 1.0}}{5\,A^2}
-3\,\frac{B_{\sssize 1.0}+A_{\sssize 0.1}+3\,B\,Q}{5\,A}
+\frac{6}{5}\,R.
\endaligned
\mytag{4.16}
$$
If both $A$ and $B$ are nonzero, then both formulas \mythetag{4.15} and
\mythetag{4.16} are applicable. The quantities $\varphi_1$ and $\varphi_2$ 
do not form a pseudotensorial field. They are used for producing the
connection components \mythetag{4.12}.\par
\head
5. Derived structures common for all cases 
of intermediate degeneration. 
\endhead
     The auxiliary structures given by the fields $M$, $N$, 
$\boldsymbol\gamma$ and the connection components $\varGamma^k_{ij}$ are
complemented by some more structures common for all cases of
{\bf intermediate degeneration}. It is known that each affine connection
produces its curvature tensor (see \mycite{13} and \mycite{14}). 
The well-known formula expressing the curvature tensor through the 
connection components is written as
$$
\hskip -2em
R^k_{rij}=
\frac{\partial\varGamma^k_{jr}}{\partial x^i}
-\frac{\partial\varGamma^k_{ir}}{\partial x^j}+
\sum^2_{q=1}\varGamma^k_{iq}\varGamma^q_{jr}-
\sum^2_{q=1}\varGamma^k_{jq}\varGamma^q_{ir}.
\mytag{5.1}
$$
The type of the field with the components \mythetag{5.1} is $(1,3)$, its
weight is zero, i\.\,e\. it is a tensorial field.\par
     The curvature tensor $R$ in \mythetag{5.1} has one upper index and three
lower indices. Contracting it with respect to the upper index and the fist
lower index, we get 
$$
\hskip -2em
\omega_{ij}=\sum^2_{k=1}R^k_{kij}.
\mytag{5.2}
$$
The quantities \mythetag{5.2} are components of a skew-symmetric tensorial
field. The quantities \mythetag{5.2} were first introduced through the quantities 
$\varphi_i$ in \mycite{3}: 
$$
\hskip -2em
\omega_{ij}=\frac{\partial\varphi_i}{\partial x^j}
-\frac{\partial\varphi_j}{\partial x^i}.
\mytag{5.3}
$$
One can verify that \mythetag{5.3} yields the same result as \mythetag{5.2}. The
quantities \mythetag{5.3} then were used in order to define a pseudoscalar field
$\Omega$ (see \myciterange{3}{3}{, }{4}):
$$
\hskip -2em
\Omega=\frac{5}{6}\sum^2_{i=1}\sum^2_{j=1}\omega_{ij}\,d^{ij}. 
\mytag{5.4}
$$
The weight of the field $\Omega$ in \mythetag{5.4} is equal to $1$. Here are
explicit formulas for $\Omega$:
$$
\align
&\hskip -2em
\aligned
\Omega&=\frac{2\,A\,B_{\sssize 0.1}(A\,S-B_{\sssize 0.1})}{B^3}
+\frac{(2\,A_{\sssize 0.1}-3\,A\,R)\,B_{\sssize 0.1}}{B^2}+\\
\vspace{1ex}
&+\frac{(B_{\sssize 1.0}-2\,A_{\sssize 0.1})\,A\,S}{B^2}
+\frac{A\,B_{\sssize 0.2}-A^2\,S_{\sssize 0.1}}{B^2}-
\frac{A_{\sssize 0.2}}{B}+\\
\vspace{1ex}
&+\frac{3\,A_{\sssize 0.1}\,R+3\,A\,R_{\sssize 0.1}
-A_{\sssize 1.0}\,S-A\,S_{\sssize 1.0}}{B}
+R_{\sssize 1.0}-2\,Q_{\sssize 0.1},
\endaligned
\mytag{5.5}\\
\vspace{1ex}
&\hskip -2em
\aligned
\Omega&=\frac{2\,B\,A_{\sssize 1.0}(B\,P+A_{\sssize 1.0})}{A^3}
-\frac{(2\,B_{\sssize 1.0}+3\,B\,Q)\,A_{\sssize 1.0}}{A^2}+\\
\vspace{1ex}
&+\frac{(A_{\sssize 0.1}-2\,B_{\sssize 1.0})\,B\,P}{A^2}
-\frac{B\,A_{\sssize 2.0}+B^2\,P_{\sssize 1.0}}{A^2}+
\frac{B_{\sssize 2.0}}{A}+\\
\vspace{1ex}
&+\frac{3\,B_{\sssize 1.0}\,Q+3\,B\,Q_{\sssize 1.0}
-B_{\sssize 0.1}\,P-B\,P_{\sssize 0.1}}{A}
+Q_{\sssize 0.1}-2\,R_{\sssize 1.0}.
\endaligned
\mytag{5.6}
\endalign
$$
The formula \mythetag{5.5} applies in the case $B\neq 0$. If $A\neq 0$, we
apply the formula \mythetag{5.6}.\par
     Note that we use the formula \mythetag{5.1} in two-dimensional case. In
two-dimensional case the curvature tensor of any connection is presented as
$$
\hskip -2em
R^k_{qij}=R^k_q\,d_{ij}\text{, \ where \ }
R^k_q=\frac{1}{2}\sum^2_{i=1}\sum^2_{j=1}R^k_{qij}\,d^{ij}.
\mytag{5.7}
$$
The type of the pseudotensorial field $R$ with the components $R^k_q$ in
\mythetag{5.7} is $(1,1)$, its weight is $1$. The field $R$ is a pseudooperator
field. There are two pseudoscalar fields associated with the field $R$ --- 
its trace $\tr(R)$ and its determinant $\det(R)$. The trace $\tr(R)$ is given 
by the following formula:
$$
\hskip -2em
\tr(R)=\frac{3}{5}\,\Omega. 
\mytag{5.8}
$$
Unlike $\tr(R)$ in \mythetag{5.8}, the determinant $\det(R)$ is a new field.
It was not studied in \myciterange{3}{3}{, }{4}, though the curvature tensor
\mythetag{5.1} and the pseudooperator field $R$ from \mythetag{5.7} were
considered in \mycite{4}. The eigenvalues of the pseudooperator field $R$ were
also considered in \mycite{4}. They are given by the following characteristic
equation:
$$
\hskip -2em
\lambda^2-\tr(R)\,\lambda+\det(R)=0. 
\mytag{5.9}
$$
\par
     Let's denote through $\goth D$ the discriminant of the quadratic equation
\mythetag{5.9}. This discriminant is calculated by means of the following formula:
$$
\hskip -2em
\goth D=\tr(R)^2-4\,\det(R)=\frac{9}{25}\,\Omega^2-4\,\det(R). 
\mytag{5.10}
$$
The determinant $\det(R)$ and the discriminant $\goth D$ in \mythetag{5.10} both
are pseudoscalar fields of the weight $2$. Both of these two pseudoscalar fields 
are defined in all cases of {\bf intermediate degeneration}.\par
\head
6. Umbilical and zero curvature equations.
\endhead
     The term {\bf umbilical points} in metric geometry of two-dimensional surfaces
in a three-dimensional Euclidean space denotes those points of a surface where two 
principal curvatures $\lambda_1$ and $\lambda_2$ are equal to each other (see
\mycite{15}). The principal curvatures are eigenvalues of a symmetric operator 
defined by the first fundamental form and the second fundamental form of a surface
(see \mycite{16}). Using this analogy, we introduce the following definition
of umbilical equations of the form \mythetag{1.1}.
\mydefinition{6.1} An umbilical equation of the form \mythetag{1.1} is an equation
of that form whose discriminant field $\goth D$ in \mythetag{5.10} is identically 
equal to zero. 
\enddefinition
     The equality $\goth D=0$ for the discriminant field \mythetag{5.10} means
that two eigenvalues $\lambda_1$ and $\lambda_2$ of the curvature pseudooperator 
field $R$ solving the characteristic equation \mythetag{5.9} are equal to each
other thus supporting the above analogy with two-dimensional surfaces.\par
     In \mycite{12} Yu\.\,Yu\.~Bagderina's classification from \mycite{11} was
compared with the previously existing classification from \myciterange{3}{3}{, }{4}. 
Among other results in \mycite{12} the following equality for Bagderina's 
pseudoscalar field $j^{\,\sssize\text{Bgd}}_5$ was derived: 
$$
\hskip -2em
j^{\,\sssize\text{Bgd}}_5=-125\,\det(R)+\frac{45}{4}\,\Omega^2.
\mytag{6.1}
$$
Comparing \mythetag{6.1} with \mythetag{5.10}, we see that 
$$
\pagebreak
\hskip -2em
j^{\,\sssize\text{Bgd}}_5=\frac{125}{4}\,\goth D. 
\mytag{6.2}
$$
The field $j^{\,\sssize\text{Bgd}}_5$ is defined by Yu\.\,Yu\.~Bagderina on
page 27 of her paper \mycite{11}. It is used in item 2 of her Theorem 2 in 
order to define one of the two basic invariants:
$$
\hskip -2em
I^{\,\sssize\text{Bgd}}_2=\frac{j^{\,\sssize\text{Bgd}}_5}
{(j^{\,\sssize\text{Bgd}}_0)^2}.
\mytag{6.3}
$$
The same field $j^{\,\sssize\text{Bgd}}_5$ is used in item 4 of 
Yu\.\,Yu\.~Bagderina's Theorem 2 in \mycite{11} again in order to define one 
of the two her basic invariants:
$$
\hskip -2em
I^{\,\sssize\text{Bgd}}_2=\frac{j^{\,\sssize\text{Bgd}}_5}
{(j^{\,\sssize\text{Bgd}}_1)^{1/2}}.
\mytag{6.4}
$$
Along with its very restrictive use in \mythetag{6.3} to \mythetag{6.4}, 
the field $j^{\,\sssize\text{Bgd}}_5$ from \mythetag{6.2} is used in the
following lemma.
\mylemma{6.1} Umbilical equations constitute that very class where 
Yu\.\,Yu\.~Bagderina's pseudoscalar field $j^{\,\sssize\text{Bgd}}_5$ is 
identically zero: $j^{\,\sssize\text{Bgd}}_5=0$.
\endproclaim
    Apart from $\lambda_1=\lambda_2$, one can write two other relationships
for principal curvatures $\lambda_1$ and $\lambda_2$, which are symmetric 
with respect to them:
$$
\xalignat 2
&\hskip -2em
\frac{\lambda_1+\lambda_2}{2}=0,
&&\lambda_1\,\lambda_2=0.
\mytag{6.5}
\endxalignat
$$
In metric geometry the expressions in the left hand sides of the formulas 
\mythetag{6.5} are called mean curvature and Gaussian curvature respectively 
(see \mycite{17} and \mycite{18}). Since
$$
\xalignat 2
&\hskip -2em
\frac{\lambda_1+\lambda_2}{2}=\frac{\tr(R)}{2},
&&\lambda_1\,\lambda_2=\det(R),
\mytag{6.6}
\endxalignat
$$
using the analogy to metric geometry, we can formulate two definitions. 
\mydefinition{6.2} A zero mean curvature equation of the form \mythetag{1.1} 
is an equation of that form whose trace field $\tr(R)$ is identically 
equal to zero. 
\enddefinition
\mydefinition{6.3} A zero Gaussian curvature equation of the form \mythetag{1.1} 
is an equation of that form whose determinant field $\det(R)$ is identically 
equal to zero. 
\enddefinition
\head
7. Zero mean curvature equations\\
in the first case of intermediate degeneration.
\endhead
    Comparing \mythetag{5.8} with \mythetag{6.6} and taking into account 
Definition~\mythedefinition{6.2}, we see that zero mean curvature equations
are given by the condition
$$
\hskip -2em
\Omega=0.
\mytag{7.1}
$$
The first case of intermediate degeneration is defined by the conditions
$$
\xalignat 3
&\hskip -2em
\boldsymbol\alpha\neq 0,
&&F=0,
&&M\neq 0
\quad
\mytag{7.2}
\endxalignat
$$
(see \myciterange{3}{3}{, }{4}). The first two of the conditions 
\mythetag{7.2} follow from \mythetag{3.11}. They are common for all cases 
of intermediate degeneration. The third condition is specific to the first 
case of intermediate degeneration. Due to \mythetag{4.11} it implies
$\boldsymbol\alpha\nparallel\boldsymbol\gamma$.\par
     Combining \mythetag{7.1} and \mythetag{7.2}, we find that zero mean
curvature equations in the first case of intermediate degeneration are given
by the following conditions:
$$
\xalignat 4
&\hskip -2em
\boldsymbol\alpha\neq 0,
&&F=0,
&&M\neq 0,
&&\Omega=0.
\quad
\mytag{7.3}
\endxalignat
$$
Now let's refer to \mycite{12} where two classifications from 
\myciterange{3}{3}{, }{4}) and \mycite{11} were compared. As a result of this 
comparison we have the following relationships
$$
\xalignat 4
&\hskip -2em
\beta^{\,\sssize\text{Bgd}}_1=\alpha_1
&&J^{\,\sssize\text{Bgd}}_0=-F^5,
&&j^{\,\sssize\text{Bgd}}_1=\frac{5}{2}\,M,
&&j^{\,\sssize\text{Bgd}}_0=-3\,\Omega\,.
\quad
\mytag{7.4}
\endxalignat
$$
Comparing \mythetag{7.3} and \mythetag{7.4} with item 4 of Theorem 2 in 
\mycite{11}, we can formulate the following comparison theorem.\par
\mytheorem{7.1} Bagderina's type four equations from \mycite{11} coincide
with the subclass of zero mean curvature equations within the first
case of intermediate degeneration in the classification from
\myciterange{3}{3}{, }{4}. 
\endproclaim
In terms of the notations introduced in the table on page 15 of \mycite{12}
we have 
$$
\hskip -2em
\text{\rm BgdET4}=\text{\rm ShrID1}\cap\text{\rm BgdET4}\subset
\text{\rm ShrID1}.
\mytag{7.5}
$$
The class {\rm ShrID1\,$\cap$\,BgdET4} in \mythetag{7.5} and in
Theorem~\mythetheorem{7.1} is complementary to the intersection class 
{\rm ShrID1\,$\cap$\,BgdET2} considered in \mycite{12}:
$$
\hskip -2em
\text{\rm ShrID1}=\bigl(\text{\rm ShrID1}\cap\text{\rm BgdET4}\bigl)
\cup\bigl(\text{\rm ShrID1}\cap\text{\rm BgdET2}\bigr).
\mytag{7.6}
$$
Due to \mythetag{7.6} the class of equations of the first case of 
intermediate degeneration {\rm ShrID1} does not intersect with the
classes other than {\rm BgdET2} and {\rm BgdET4} in Bagderina's
classification from \mycite{11}.\par
     Relying on Theorem~\mythetheorem{7.1}, below we continue studying the
intersection class {\rm ShrID1\,$\cap$\,BgdET4}. In \mycite{11} 
Yu\.\,Yu\.\,Bagderina defines two invariant differentiation operators. 
The first of them is given by the formula 
$$
\hskip -2em
\Cal D^{\,\sssize\text{Bgd}}_1=\frac{\beta^{\,\sssize\text{Bgd}}_2}
{(\mu^{\,\sssize\text{Bgd}}_1)^2}
\,\frac{\partial}{\partial x}
-\frac{\beta^{\,\sssize\text{Bgd}}_1}
{(\mu^{\,\sssize\text{Bgd}}_1)^2}\,\frac{\partial}{\partial y}.
\mytag{7.7}
$$
The quantity $\mu^{\,\sssize\text{Bgd}}_1$ in \mythetag{7.7} is given
by one of the formulas \thetag{2.12} from \mycite{11}:
$$
\hskip -2em
\mu^{\,\sssize\text{Bgd}}_1=\bigl(j^{\,\sssize\text{Bgd}}_1\bigr)^{1/4}.
\mytag{7.8}
$$
Comparing the formula \mythetag{7.8} with  \mythetag{7.4}, we can write 
$$
\hskip -2em
\mu^{\,\sssize\text{Bgd}}_1=\root{4}\of{\frac{5\,M}{2}}.
\mytag{7.9}
$$
The quantities $\beta^{\,\sssize\text{Bgd}}_1$ and $\beta^{\,\sssize\text{Bgd}}_2$
in \mythetag{7.7} coincide with the components \mythetag{3.1} of the pseudoscalar
field $\boldsymbol\alpha$ (see Lemma~3.2 in \mycite{19} or Lemma~5.1 in 
\mycite{12}): 
$$
\pagebreak
\xalignat 2
&\hskip -2em
\beta^{\,\sssize\text{Bgd}}_1=\alpha_1=A,
&&\beta^{\,\sssize\text{Bgd}}_2=\alpha_2=B.
\mytag{7.10}
\endxalignat
$$
Raising indices in \mythetag{7.10} according to \mythetag{3.5} and then applying 
\mythetag{7.10} and \mythetag{7.9} to \mythetag{7.7}, we derive the following
formula for $\Cal D^{\,\sssize\text{Bgd}}_1$:
$$
\hskip -2em
\Cal D^{\,\sssize\text{Bgd}}_1=\sqrt{\frac{2}{5\,M}}\biggl(\alpha^1
\,\frac{\partial}{\partial x}+\alpha^2\,\frac{\partial}{\partial y}\biggr).
\mytag{7.11}
$$
\mylemma{7.1} If the conditions \mythetag{7.3} are fulfilled, i\.\,e\.
within the intersection class\/ {\rm ShrID1\,$\cap$\,BgdET4} coinciding with
the class of Bagderina's type four equations, Bagderina's invariant differentiation 
operator $\Cal D^{\,\sssize\text{Bgd}}_1$ is expressed through the pseudovectorial
field $\boldsymbol\alpha$ and pseudoscalar field $M$ from \mycite{3} by means of 
the formula \mythetag{7.11}.
\endproclaim
    The use of non-integer power exponents with even denominators by Bagderina 
in \mythetag{6.4} and \mythetag{7.8} is a bad practice. In order to avoid such 
a practice in \myciterange{3}{3}{, }{4} covariant differentiation operators 
were considered instead of invariant differentiations. They are produced
by pseudovectorial fields according to the following patterns:
$$
\xalignat 2
&\hskip -2em
\nabla_{\boldsymbol\alpha}=\alpha^1\,\nabla_1+\alpha^2\,\nabla_2,
&&\nabla_{\boldsymbol\gamma}=\gamma^1\,\nabla_1+\gamma^2\,\nabla_2
\mytag{7.12}
\endxalignat
$$
(see \thetag{6.13} in \mycite{3} and \thetag{5.2} in \mycite{4}). The covariant 
derivatives in \mythetag{7.12} extend partial derivatives from \mythetag{7.7}.
They are defined by means of the formula
$$
\hskip -2em
\aligned
\nabla_kF^{i_1\ldots\,i_r}_{j_1\ldots\,j_s}&=
\frac{\partial F^{i_1\ldots\,i_r}_{j_1\ldots\,j_s}}{\partial x^k}+
\sum^r_{n=1}\sum^2_{v_n=1}\varGamma^{i_n}_{k\,v_n}\,
F^{i_1\ldots\,v_n\ldots\,i_r}_{j_1\ldots\,j_s}-\\
&-\sum^s_{n=1}\sum^2_{w_n=1}\varGamma^{w_n}_{k\,j_n}\,
F^{i_1\ldots\,i_r}_{j_1\ldots\,w_n\ldots\,j_s}+
m\,\varphi_k\,F^{i_1\ldots\,i_r}_{j_1\ldots\,j_s}
\endaligned
\mytag{7.13}
$$
(see \thetag{6.11} in \mycite{3} or \thetag{4.23} in \mycite{4}). The 
covariant derivative $\nabla_k$ in \mythetag{7.13} is applied to a 
pseudotensorial field of the type $(r,s)$ and the weight $m$. The connection 
components $\varGamma^k_{ij}$ in \mythetag{7.13} are defined by \mythetag{4.12}. 
They are canonically associated with a given equation \mythetag{1.1} in all
cases of intermediate degeneration.\par
     Covariant differentiation operators like \mythetag{7.12} are preferable 
with respect to invariant differentiation operators like \mythetag{7.7}. They can 
be applied to any pseudotensorial invariants, not only to scalar ones. Fortunately
invariant differentiation operators can be extended to covariant differentiations.
In the case of \mythetag{7.11} we have
$$
\hskip -2em
\Cal D^{\,\sssize\text{Bgd}}_1=\sqrt{\frac{2}{5\,M}}\,\bigl(\alpha^1
\,\nabla_1+\alpha^2\,\nabla_2\bigr)=\sqrt{\frac{2}{5\,M}}\,
\nabla_{\boldsymbol\alpha}.
\mytag{7.14}
$$
Due to \mythetag{7.14} we can reformulate Lemma~\mythelemma{7.1} as follows.
\mylemma{7.2} If the conditions \mythetag{7.3} are fulfilled, i\.\,e\.
within the intersection class\/ {\rm ShrID1\,$\cap$\,BgdET4} coinciding with
the class of Bagderina's type four equations, Bagderina's invariant differentiation 
operator $\Cal D^{\,\sssize\text{Bgd}}_1$ is expressed through the pseudovectorial
field $\boldsymbol\alpha$ and pseudoscalar field $M$ from \mycite{3} by means of 
the formula \mythetag{7.14}.
\endproclaim
     The second invariant differentiation operator $\Cal D^{\,\sssize\text{Bgd}}_1$
introduced by Yu\.\,Yu\.\,Bagde\-rina in \mycite{11} upon replacing partial derivatives
by covariant derivatives is 
$$
\hskip -2em
\Cal D^{\,\sssize\text{Bgd}}_2=
\biggl(\mu^{\,\sssize\text{Bgd}}_2\,\beta^{\,\sssize\text{Bgd}}_2
-3\,\frac{\mu^{\,\sssize\text{Bgd}}_1}{\beta^{\,\sssize\text{Bgd}}_1}
\biggr)\,\nabla_1
-\mu^{\,\sssize\text{Bgd}}_2\,\beta^{\,\sssize\text{Bgd}}_1\,\nabla_2.
\mytag{7.15}
$$
Note that the formula \thetag{6.14} for $\Cal D^{\,\sssize\text{Bgd}}_2$ in 
\mycite{12} is mistyped. The author apologizes for this typo and presents the
correct formula \mythetag{7.15}.\par 
    The quantity $\mu^{\,\sssize\text{Bgd}}_2$ in \mythetag{7.15} is given by 
one of the formulas \thetag{2.12} in \mycite{11}:
$$
\mu^{\,\sssize\text{Bgd}}_2=\frac{5\,j^{\,\sssize\text{Bgd}}_2}
{2\,\bigl(j^{\,\sssize\text{Bgd}}_1\bigr)^{3/4}}
\mytag{7.16}
$$
The quantity $j^{\,\sssize\text{Bgd}}_2$ in \mythetag{7.16} is expressed 
through the coefficients of the equation \mythetag{1.1} by means of a series of
auxiliary notations (see \mycite{11}). We reproduce them here without exceptions
and reductions for the sake of completeness:
$$
\gather
\aligned
&\aligned
j^{\,\sssize\text{Bgd}}_2&=\frac{1}{\beta^{\,\sssize\text{Bgd}}_1}
\biggl(\delta^{\,\sssize\text{Bgd}}_{20}
-\frac{\beta^{\,\sssize\text{Bgd}}_2}{\beta^{\,\sssize\text{Bgd}}_1}
\,\delta^{\,\sssize\text{Bgd}}_{10}\biggr)\,+\\
&+\,\frac{\gamma^{\sssize\text{Bgd}}_{10}}{5\,(\beta^{\,\sssize\text{Bgd}}_1)^2}
\biggl(7\,\frac{\beta^{\,\sssize\text{Bgd}}_2}{\beta^{\,\sssize\text{Bgd}}_1}
\,\gamma^{\sssize\text{Bgd}}_{10}-6\,\gamma^{\sssize\text{Bgd}}_{20}
-\gamma^{\sssize\text{Bgd}}_{11}\biggr),
\endaligned\\
\vspace{1ex}
&j^{\,\sssize\text{Bgd}}_3=\frac{3}{5}\biggl(
\frac{\delta^{\,\sssize\text{Bgd}}_{10}}{(\beta^{\,\sssize\text{Bgd}}_1)^3}
-\frac{6\,(\gamma^{\sssize\text{Bgd}}_{10})^2}
{5\,(\beta^{\,\sssize\text{Bgd}}_1)^4}\biggr),
\endaligned
\mytag{7.17}\\
\vspace{2ex}
\hskip -2em
\aligned
&\delta^{\,\sssize\text{Bgd}}_{10}=\partial_x\gamma^{\sssize\text{Bgd}}_{10}
-2\,Q\,\gamma^{\sssize\text{Bgd}}_{10}+P\,(\gamma^{\sssize\text{Bgd}}_{20}
+\gamma^{\sssize\text{Bgd}}_{11})-5\,\alpha^{\sssize\text{Bgd}}_0\,
\beta^{\,\sssize\text{Bgd}}_1,\\
\vspace{1ex}
&\delta^{\,\sssize\text{Bgd}}_{20}=\partial_x\gamma^{\sssize\text{Bgd}}_{20}
-R\,\gamma^{\sssize\text{Bgd}}_{10}+P\,\gamma^{\sssize\text{Bgd}}_{21}
-4\,\alpha^{\sssize\text{Bgd}}_1\,\beta^{\,\sssize\text{Bgd}}_1
-\alpha^{\sssize\text{Bgd}}_0\,\beta^{\,\sssize\text{Bgd}}_2,\\
\endaligned
\mytag{7.18}\\
\vspace{2ex}
\hskip -2em
\aligned
&\gamma^{\sssize\text{Bgd}}_{10}=\partial_x\beta^{\,\sssize\text{Bgd}}_1
-Q\,\beta^{\,\sssize\text{Bgd}}_1+P\,\beta^{\,\sssize\text{Bgd}}_2,\\
\vspace{1ex}
&\gamma^{\sssize\text{Bgd}}_{11}=\partial_x\beta^{\,\sssize\text{Bgd}}_2
-R\,\beta^{\,\sssize\text{Bgd}}_1+Q\,\beta^{\,\sssize\text{Bgd}}_2,\\
\vspace{1ex}
&\gamma^{\sssize\text{Bgd}}_{20}=\partial_y\beta^{\,\sssize\text{Bgd}}_1
-R\,\beta^{\,\sssize\text{Bgd}}_1+Q\,\beta^{\,\sssize\text{Bgd}}_2,\\
\vspace{1ex}
&\gamma^{\sssize\text{Bgd}}_{21}=\partial_y\beta^{\,\sssize\text{Bgd}}_2
-S\,\beta^{\,\sssize\text{Bgd}}_1+R\,\beta^{\,\sssize\text{Bgd}}_2,
\endaligned
\mytag{7.19}\\
\vspace{2ex}
\hskip -2em
\aligned
&\alpha^{\sssize\text{Bgd}}_0=Q_{\sssize 1.0}-P_{\sssize 0.1}+2\,P\,R-2\,Q^2,\\
\vspace{1ex}
&\alpha^{\sssize\text{Bgd}}_1=R_{\sssize 1.0}-Q_{\sssize 0.1}+P\,S-Q\,R,\\
\vspace{1ex}
&\alpha^{\sssize\text{Bgd}}_2=S_{\sssize 1.0}-R_{\sssize 0.1}+2\,Q\,S-2\,R^2.
\endaligned
\mytag{7.20}
\endgather
$$
We do no provide Bagderina's expressions for $j^{\,\sssize\text{Bgd}}_1$,
$\beta^{\,\sssize\text{Bgd}}_1$, $\beta^{\,\sssize\text{Bgd}}_2$ since they
can be calculated using \mythetag{7.4} and \mythetag{7.10}. Using the formulas
\mythetag{7.16}, \mythetag{7.17}, \mythetag{7.18}, \mythetag{7.19}, 
\mythetag{7.20} and applying them to \mythetag{7.15}, one can derive the
following formula for $\Cal D^{\,\sssize\text{Bgd}}_2$:
$$
\hskip -2em
\Cal D^{\,\sssize\text{Bgd}}_2=-3\,\root{4}\of{\frac{5}{2\,M^3}}
\,\bigl(\gamma^1\,\nabla_1+\gamma^2\,\nabla_2\bigr)=
-3\,\root{4}\of{\frac{5}{2\,M^3}}\,\nabla_{\boldsymbol\gamma}.
\mytag{7.21}
$$
\mylemma{7.3} If the conditions \mythetag{7.3} are fulfilled, i\.\,e\.
within the intersection class\/ {\rm ShrID1\,$\cap$\,BgdET4} coinciding with
the class of Bagderina's type four equations, Bagderina's invariant differentiation 
operator $\Cal D^{\,\sssize\text{Bgd}}_2$ is expressed through the pseudovectorial
field $\boldsymbol\gamma$ and pseudoscalar field $M$ from \mycite{3} by means of 
the formula \mythetag{7.21}.
\endproclaim
     Lemmas~\mythelemma{7.1}, \mythelemma{7.2}, and \mythelemma{7.3} are proved 
by direct calculations using some symbolic algebra package. In my case that was Maple\footnotemark.\par 
\footnotetext{\ Maple is a trademark of Waterloo Maple Inc.}
\head
8. Scalar invariants.
\endhead
     The condition $\Omega=0$ in \mythetag{7.3} does not specify a special subcase 
within the first case of intermediate degeneration. The subcase $\Omega=0$ is treated 
regularly, though one of the three invariants $I_1$, $I_2$, $I_3$ from 
\myciterange{3}{3}{, }{4}) does vanish:
$$
\hskip -2em
I_2=\frac{\Omega^2}{N}=0.
\mytag{8.1}
$$
The invariant $I_1$ does not vanish due to \mythetag{7.3}. It is given by the formula
$$
\hskip -2em
I_1=\frac{M}{N^2}.
\mytag{8.2}
$$
The third invariant $I_3$ is introduced through the coefficients in the expansions 
$$
\xalignat 2
&\hskip -2em
\nabla_{\boldsymbol\alpha}\boldsymbol\alpha=\Gamma^1_{11}\,\boldsymbol\alpha+
\Gamma^2_{11}\,\boldsymbol\gamma,
&&\nabla_{\boldsymbol\alpha}\boldsymbol\gamma=\Gamma^1_{12}\,\boldsymbol\alpha+
\Gamma^2_{12}\,\boldsymbol\gamma,\\
\vspace{-1.7ex}
&&&\mytag{8.3}\\
\vspace{-1.7ex}
&\hskip -2em
\nabla_{\boldsymbol\gamma}\boldsymbol\alpha=\Gamma^1_{21}\,\boldsymbol\alpha+
\Gamma^2_{21}\,\boldsymbol\gamma,
&&\nabla_{\boldsymbol\gamma}\boldsymbol\gamma=\Gamma^1_{22}\,\boldsymbol\alpha+
\Gamma^2_{22}\,\boldsymbol\gamma
\endxalignat
$$
(see \thetag{6.13} in \mycite{3} or \thetag{5.2} in \mycite{4}). As it was shown 
in \mycite{3}, only one of the coefficients $\Gamma^1_{11}$, $\Gamma^2_{11}$, 
$\Gamma^1_{12}$, $\Gamma^2_{12}$, $\Gamma^1_{21}$, $\Gamma^2_{21}$, 
$\Gamma^1_{22}$, $\Gamma^2_{22}$ does matter. This is the coefficient 
$\Gamma^1_{22}$. It defines the invariant $I_3$ by means of the formula
$$
\hskip -2em
I_3=\frac{\Gamma^1_{22}\,N^2}{M^2}
\mytag{8.4}
$$
(see \thetag{6.20} in \mycite{3}). The coefficients $\Gamma^2_{11}$ and
$\Gamma^1_{21}$ in \mythetag{8.3} are equal to zero
$$
\xalignat 2
&\hskip -2em
\Gamma^2_{11}=0,
&&\Gamma^1_{21}=0.
\mytag{8.5}
\endxalignat
$$ 
The coefficients $\Gamma^1_{11}$ and $\Gamma^2_{21}$ are expressed through 
the pseudoscalar field $N$:
$$
\xalignat 2
&\hskip -2em
\Gamma^1_{11}=-\frac{3}{5}\,N,
&&\Gamma^2_{21}=-\frac{3}{5}\,N.
\mytag{8.6}
\endxalignat
$$ 
The coefficients $\Gamma^2_{12}$, $\Gamma^2_{22}$, and $\Gamma^1_{12}$
obey more complicated relationships
$$
\gather
\hskip -2em
I_1\,\Gamma^2_{12}=I_4\,N+\frac{3}{5}\,I_1\,N+2\,I_1^2\,N,
\mytag{8.7}\\
\vspace{2ex}
\hskip -2em
\aligned
 \bigl(I_1\,\Gamma^2_{22}\bigr)^4&+\bigl(I_7\,N^3\bigr)^2
 +\bigl(16\,I_2\,N^3\,{I_1}^4\bigr)^2=\\
 \vspace{1ex}
 &=32\,I_7\,N^6\,I_2\,{I_1}^4
 +2\,\bigl(I_7\,N^3+16\,I_2\,N^3\,{I_1}^4\bigr)\,
 \bigl(I_1\,\Gamma^2_{22}\bigr)^2,
\endaligned
\mytag{8.8}\\
\vspace{2ex}
\hskip -2em
\Gamma^1_{12}=-\Gamma^2_{22}.
\mytag{8.9}
\endgather
$$
(see \thetag{6.22}, \thetag{6.23}, \thetag{6.15} and  in \mycite{3} or 
\thetag{5.6} and \thetag{5.7} in \mycite{4}). Unfortunately the formula 
\thetag{6.22} in \mycite{3} is mistyped and then copied to \thetag{5.6} 
in \mycite{4}. The minus signs in the right hand side of this formula should 
be altered for pluses. The formula \mythetag{8.7} here is a corrected 
version of this formula.\par
     The quantities $I_4$ and $I_7$ in \mythetag{8.7} and \mythetag{8.8} are
higher order invariants. They are taken from the following formulas:
$$
\xalignat 3
&\hskip -2em
I_4=\frac{\nabla_{\boldsymbol\alpha}I_1}{N},
&&I_5=\frac{\nabla_{\boldsymbol\alpha}I_2}{N},
&&I_6=\frac{\nabla_{\boldsymbol\alpha}I_3}{N},\\
\vspace{-1.5ex}
\mytag{8.10}\\
\vspace{-1.5ex}
&\hskip -2em
I_7=\frac{(\nabla_{\boldsymbol\gamma}I_1)^2}{N^3},
&&I_8=\frac{(\nabla_{\boldsymbol\gamma}I_2)^2}{N^3},
&&I_9=\frac{(\nabla_{\boldsymbol\gamma}I_3)^2}{N^3}.
\endxalignat
$$
More higher order invariants can be produced recursively 
$$
\xalignat 2
&\hskip -2em
I_{k+3}=\frac{\nabla_{\boldsymbol\alpha}I_k}{N},
&&I_{k+6}=\frac{\bigl(\nabla_{\boldsymbol\gamma}I_k\bigr)^2}{N^3}.
\mytag{8.11}
\endxalignat
$$
The recurrent formulas \mythetag{8.11} should be applied in triples in some 
commonly negotiated order (see \thetag{6.21} in \mycite{3}). As a result we 
shall have an infinite series of scalar invariants associated with the 
equation \mythetag{1.1}.\par
     In \mycite{11} Yu\.~Yu\. Bagderina introduces her own basic invariants
for the class {\rm BgdET4} coinciding with the intersection class 
{\rm ShrID1\,$\cap$\,BgdET4} and given by the conditions \mythetag{7.3}. 
Her first invariant $I^{\,\sssize\text{Bgd}}_1$ is given by the formula
$$
\hskip -2em
I^{\,\sssize\text{Bgd}}_1=\frac{\Gamma^{\,\sssize\text{Bgd}}_0}
{\beta^{\,\sssize\text{Bgd}}_1\,(j^{\,\sssize\text{Bgd}}_1)^{1/2}}.
\mytag{8.12}
$$
The fields $\beta^{\,\sssize\text{Bgd}}_1$ and $j^{\,\sssize\text{Bgd}}_1$
in \mythetag{8.12} are taken from \mythetag{7.10} and \mythetag{7.4}. As for the
field $\Gamma^{\,\sssize\text{Bgd}}_0$, in \mycite{19} it was shown that the
$\Gamma^{\,\sssize\text{Bgd}}_0$ coincides with the first component of the
pseudocovectorial field $\boldsymbol\beta$ in \mythetag{3.2} (see Lemma~3.5
in \mycite{19} or Lemma~5.2 in \mycite{12}). Therefore from \mythetag{4.1}
or from \mythetag{4.2} in our present case we derive
$$
\hskip -2em
\frac{\Gamma^{\,\sssize\text{Bgd}}_0}
{\beta^{\,\sssize\text{Bgd}}_1}=3\,N.
\mytag{8.13}
$$
Applying \mythetag{8.13}, \mythetag{7.4} and \mythetag{8.2} to \mythetag{8.12}, 
we derive the following formula:
$$
\hskip -2em
I^{\,\sssize\text{Bgd}}_1=\sqrt{\frac{18\,N^2}{5\,M}}
=\sqrt{\frac{18}{5\,I_1}}.
\mytag{8.14}
$$
The following lemma is formulated for further references. 
\mylemma{8.1} If the conditions \mythetag{7.3} are fulfilled, i\.\,e\.
within the intersection class\/ {\rm ShrID1\,$\cap$\,BgdET4} coinciding with
the class of Bagderina's type four equations, Bagderina's basic invariant
$I^{\,\sssize\text{Bgd}}_1$ is expressed through the invariant $I_1$
from \mycite{3} by means of the formula \mythetag{8.14}.
\endproclaim
     The formula in the right hand side of \mythetag{8.14} is in agreement 
with the first formula \thetag{7.2} presented by Yu\.~Yu\. Bagderina in 
\mycite{20}, where she compares her results with the previously known results 
from \myciterange{3}{3}{--}{5}. The formula \mythetag{8.1} is evidently
in agreement with the second formula \thetag{7.2} from \mycite{20}. 
The third formula \thetag{7.2} in \mycite{20} is more complicated. We
shall consider it below.\par
     Apart from \thetag{7.2}, there are the following comparison formulas
in \mycite{20}: 
$$
\xalignat 5
&A=\beta^{\,\sssize\text{Bgd}}_1,
&&B=\beta^{\,\sssize\text{Bgd}}_2,
&&F^5=-J^{\,\sssize\text{Bgd}}_0,
&&G=\Gamma^{\,\sssize\text{Bgd}}_1,
&&H=-\Gamma^{\,\sssize\text{Bgd}}_0.
\qquad\quad
\mytag{8.15}
\endxalignat
$$
The formulas \mythetag{8.15} are in agreement with the results from
\mycite{19} and \mycite{12} (see Lemma~3.2, Lemma~3.3, and Lemma~3.5
in \mycite{19}). The formulas 
$$
\xalignat 3
&\hskip -2em
N=\frac{\Gamma^{\,\sssize\text{Bgd}}_0}{3\,\beta^{\,\sssize\text{Bgd}}_1},
&&\Omega=-\frac{j^{\,\sssize\text{Bgd}}_0}{3},
&&M=\frac{2}{5}\,j^{\,\sssize\text{Bgd}}_1
\quad
\mytag{8.16}
\endxalignat
$$
are also presented in \mycite{20}. The formulas \mythetag{8.16} are in agreement 
with the formulas \mythetag{8.13} and \mythetag{7.4} in the present paper. The 
formulas \thetag{7.1} from \mycite{20} look like 
$$
\xalignat 4
&\hskip -2em
F=0
&&A\neq 0, 
&&\Omega=0, 
&&M\neq 0. 
\quad
\mytag{8.17}
\endxalignat
$$
The formulas \mythetag{8.17} are equivalent to the condition \mythetag{7.3}. 
They define the class {\rm BgdET4} coinciding with the intersection class 
{\rm ShrID1\,$\cap$\,BgdET4}.\par
     Immediately after the formulas \thetag{7.1} in \mycite{20} we see the
formulas 
$$
\xalignat 2
&\hskip -2em
\gamma^1=\frac{2\,j^{\,\sssize\text{Bgd}}_1}
{5\,\beta^{\,\sssize\text{Bgd}}_1}-\frac{\beta^{\,\sssize\text{Bgd}}_2\,
j^{\,\sssize\text{Bgd}}_2}{3},
&&\gamma^2=\frac{\beta^{\,\sssize\text{Bgd}}_1\,
j^{\,\sssize\text{Bgd}}_2}{3}.
\mytag{8.18}
\endxalignat
$$
The formulas \mythetag{8.18} are valid. They are verified by direct
calculations.\par
     Now we can proceed to the third formula \thetag{7.2} in \mycite{20}. 
It is written as follows:
$$
\hskip -2em
I_3=\frac{I^{\,\sssize\text{Bgd}}_{12}}{18}
+\frac{I^{\,\sssize\text{Bgd}}_2}{90\,I^{\,\sssize\text{Bgd}}_1}
\Bigl(5\,I^{\,\sssize\text{Bgd}}_{11}
-3\,\bigl(I^{\,\sssize\text{Bgd}}_1\bigr)^2-6\Bigr)+\frac{5}{3}.
\mytag{8.19}
$$
The invariants $I^{\,\sssize\text{Bgd}}_{12}$ and $I^{\,\sssize\text{Bgd}}_{11}$
in \mythetag{8.19} are calculated by applying the invariant differentiation
operator $\Cal D^{\,\sssize\text{Bgd}}_1$ from \mycite{11} to 
$I^{\,\sssize\text{Bgd}}_2$ and $I^{\,\sssize\text{Bgd}}_1$ respectively:
$$
\xalignat 2
&\hskip -2em
I^{\,\sssize\text{Bgd}}_{12}
=\Cal D^{\,\sssize\text{Bgd}}_1(I^{\,\sssize\text{Bgd}}_2),
&&I^{\,\sssize\text{Bgd}}_{11}
=\Cal D^{\,\sssize\text{Bgd}}_1(I^{\,\sssize\text{Bgd}}_1).
\mytag{8.20}
\endxalignat
$$
Using Lemma~\mythelemma{7.2}, the formulas \mythetag{8.20} can be rewritten 
as
$$
\xalignat 2
&\hskip -2em
I^{\,\sssize\text{Bgd}}_{12}
=\sqrt{\frac{2}{5\,M}}\,\nabla_{\boldsymbol\alpha}I^{\,\sssize\text{Bgd}}_2,
&&I^{\,\sssize\text{Bgd}}_{11}
=\sqrt{\frac{2}{5\,M}}\,\nabla_{\boldsymbol\alpha}I^{\,\sssize\text{Bgd}}_1.
\mytag{8.21}
\endxalignat
$$
The invariants $I^{\,\sssize\text{Bgd}}_2$ and $I^{\,\sssize\text{Bgd}}_1$
in \mythetag{8.20} and \mythetag{8.21} are two basic invariants defined by
Yu\.~Yu\. Bagderina in \mycite{11} in item 4 of her Theorem 2. The invariant
$I^{\,\sssize\text{Bgd}}_2$ is given by the formula \mythetag{6.4}, the 
invariant $I^{\,\sssize\text{Bgd}}_1$ is given by the formula \mythetag{8.12}.
We can apply the formula \mythetag{6.1} in order to derive an explicit formula 
for $j^{\,\sssize\text{Bgd}}_5$ in \mythetag{6.4}. The formulas \mythetag{6.1},
\mythetag{6.4}, \mythetag{8.14}, and \mythetag{8.21} are sufficient in order
to verify the formula \mythetag{8.19} by means of direct computations. It turns
out that the formula \mythetag{8.19} is written for the case where $I_3$ is
redefined as
$$
\hskip -2em
I_3\to \frac{I_3}{I_1}=\frac{\Gamma^1_{22}}{M}.
\mytag{8.22}
$$
Probably \mythetag{8.22} would be a better choice for $I_3$. But, historically
in \myciterange{3}{3}{, }{4} it was introduced in its present form \mythetag{8.4}. 
We need to rewrite the formula \mythetag{8.19} as 
$$
\hskip -2em
\frac{I_3}{I_1}=\frac{I^{\,\sssize\text{Bgd}}_{12}}{18}
+\frac{I^{\,\sssize\text{Bgd}}_2}{90\,I^{\,\sssize\text{Bgd}}_1}
\Bigl(5\,I^{\,\sssize\text{Bgd}}_{11}
-3\,\bigl(I^{\,\sssize\text{Bgd}}_1\bigr)^2-6\Bigr)+\frac{5}{3}.
\mytag{8.23}	
$$
\par
     The formula \mythetag{8.23} is valid. It is verified by means of direct
computations. Looking at the formula \mythetag{8.23} and at few other formulas
in \mycite{20} expressing $I_4$, $I_6$, $I_7$, and $I_{10}$ through her 
invariants, Yu\.\,Yu\.\,Bagderina detected that they do not comprise her basic
invariant $I^{\,\sssize\text{Bgd}}_2$, but only some derivatives of 
$I^{\,\sssize\text{Bgd}}_2$. As a result she issued a criticism saying that 
the invariants from \myciterange{3}{3}{, }{4} are impractical for the equivalence 
problem. However, she omitted the invariant $I_{9}$ in the sequence 
$I_4$, $I_6$, $I_7$, and $I_{10}$. The invariants $I_5$ and $I_8$, which are 
also omitted, are zero due to $I_2=0$ (see \mythetag{8.10} and \mythetag{8.1}).
But the invariant $I_9$ is nonzero. If Yu\.\,Yu\.\,Bagderina would not omit 
this invariant, she would have the following formula:
$$
\gathered
\sqrt{I_9}=-\frac{\sqrt{3}}{45}\,\frac{I^{\,\sssize\text{Bgd}}_{212}}
{\bigl(I^{\,\sssize\text{Bgd}}_1\bigr)^{3/2}}
-\frac{\sqrt{3}}{225}\,\frac{\bigl(5\,I^{\,\sssize\text{Bgd}}_{11}
-3\,\bigl(I^{\,\sssize\text{Bgd}}_1\bigr)^2-6\bigr)
\,I^{\,\sssize\text{Bgd}}_{22}}
{\bigl(I^{\,\sssize\text{Bgd}}_1\bigr)^{5/2}}\,+\\
+\biggl(\frac{\sqrt{3}}{450}
\,\frac{\bigl(15\,I^{\,\sssize\text{Bgd}}_{11}
+10\,\bigl(I^{\,\sssize\text{Bgd}}_1\bigr)^2-6\bigr)
\,I^{\,\sssize\text{Bgd}}_{21}}
{\bigl(I^{\,\sssize\text{Bgd}}_1\bigr)^{7/2}}
-\frac{\sqrt{3}}{45}
\,\frac{I^{\,\sssize\text{Bgd}}_{121}}
{\bigl(I^{\,\sssize\text{Bgd}}_1\bigr)^{5/2}}\biggr)
\,I^{\,\sssize\text{Bgd}}_2.
\endgathered
\quad
\mytag{8.24}
$$
The invariant $I^{\,\sssize\text{Bgd}}_2$ is explicitly present in
the formula \mythetag{8.24}. Opposing Yu\.\,Yu\.\,Bag\-derina, below we 
prove that, in spite of $\Omega=0$, in spite of \mythetag{8.23}, and in 
spite of other her formulas in \mycite{20}, her basic invariants 
$I^{\,\sssize\text{Bgd}}_1$ and $I^{\,\sssize\text{Bgd}}_2$ can be 
expressed through the invariants $I_1$ and $I_3$ from
\myciterange{3}{3}{, }{4} in her class {\rm BgdET4} coinciding with the 
intersection class {\rm ShrID1\,$\cap$\,BgdET4}.\par
     Let's recall that in \mycite{12} the following theorem was proved
for Bagderina's pseudoscalar field $j^{\,\sssize\text{Bgd}}_5$ from
her paper \mycite{11}. 
\mytheorem{8.1} Within the intersection class\/ {\rm ShrID1\,$\cap$\,BgdET2}, 
i\.\,e\. if the conditions $F=0$, $\boldsymbol\alpha\neq 0$,
$M\neq 0$, $\Omega\neq 0$ are fulfilled, Bagderina's pseudoscalar
field $j^{\,\sssize\text{Bgd}}_5$ is expressed through $M$, through $\Omega$, 
through the pseudoscalar fields $\Gamma^1_{11}$, $\Gamma^2_{11}$, 
$\Gamma^1_{12}$, $\Gamma^2_{12}$, $\Gamma^1_{21}$, $\Gamma^2_{21}$, 
$\Gamma^1_{22}$, $\Gamma^2_{22}$ from \mythetag{8.3}, and through covariant 
derivatives of them along the pseudovectorial fields $\boldsymbol\alpha$ and $\boldsymbol\gamma$. 
\endproclaim
     The field $j^{\,\sssize\text{Bgd}}_5$ is defined for all cases of 
intermediate degeneration. In is introduced in \mycite{11} by means of 
the following formula:
$$
\hskip -2em
j^{\,\sssize\text{Bgd}}_5=5\,\bigl(2\,j^{\,\sssize\text{Bgd}}_1\,
j^{\,\sssize\text{Bgd}}_3+(j^{\,\sssize\text{Bgd}}_2
-j^{\,\sssize\text{Bgd}}_0/6)^2\bigr).
\mytag{8.25}
$$
As it was shown in \mycite{12}, the formula \mythetag{8.25} is equivalent
to the formula \mythetag{6.1} (see \thetag{7.5} in \mycite{12}). In our
present case $j^{\,\sssize\text{Bgd}}_0=-3\,\Omega=0$ (see \mythetag{7.3}
and \mythetag{7.4}). Therefore the formula \mythetag{8.25} reduces to
the following formula:
$$
\hskip -2em
j^{\,\sssize\text{Bgd}}_5=5\,\bigl(2\,j^{\,\sssize\text{Bgd}}_1\,
j^{\,\sssize\text{Bgd}}_3+(j^{\,\sssize\text{Bgd}}_2
)^2\bigr).
\mytag{8.26}
$$
The formula \mythetag{8.26} is given in item 4 of Bagderina's Theorem 2 in
\mycite{11}. This formula defines $j^{\,\sssize\text{Bgd}}_5$ in our present 
class {\rm BgdET4} coinciding with the intersection class 
{\rm ShrID1\,$\cap$\,BgdET4}. It is equivalent to the reduced formula
\mythetag{6.1}:
$$
\hskip -2em
j^{\,\sssize\text{Bgd}}_5=-125\,\det(R). 
\mytag{8.27}
$$
The formula \mythetag{8.27} is produced from \mythetag{6.1} by setting
$\Omega=0$ in it.\par 
     Looking through the proof of Theorem~\mythetheorem{8.1} in \mycite{12},
one can see that it does not depend on $\Omega$ otherwise than through the
entry of $\Omega^2$ in \mythetag{6.1}. \pagebreak Therefore, repeating 
the arguments from \mycite{12}, we can prove the following theorem for 
Bagderina's pseudoscalar field $j^{\,\sssize\text{Bgd}}_5$ in 
\mythetag{8.26}. 
\mytheorem{8.2} Within the intersection class\/ {\rm ShrID1\,$\cap$\,BgdET4}
coinciding with {\rm BgdET4}, i\.\,e\. if the conditions 
\mythetag{7.3} are fulfilled, Bagderina's pseudoscalar field 
$j^{\,\sssize\text{Bgd}}_5$ from \mythetag{8.26} is expressed through $M$, 
through the pseudoscalar fields $\Gamma^1_{11}$, $\Gamma^2_{11}$, 
$\Gamma^1_{12}$, $\Gamma^2_{12}$, $\Gamma^1_{21}$, $\Gamma^2_{21}$, 
$\Gamma^1_{22}$, $\Gamma^2_{22}$ from \mythetag{8.3}, and through covariant 
derivatives of them along the pseudovectorial fields $\boldsymbol\alpha$ and $\boldsymbol\gamma$. 
\endproclaim
     Bagderina's invariant $I^{\,\sssize\text{Bgd}}_2$ is given by the formula
\mythetag{6.4}. The field $j^{\,\sssize\text{Bgd}}_5$ is in the numerator
of this formula, while its denominator is defined by the field $M$ due to 
\mythetag{7.4} or \mythetag{8.17}. As for the pseudoscalar fields $\Gamma^1_{11}$, 
$\Gamma^2_{11}$, $\Gamma^1_{12}$, $\Gamma^2_{12}$, $\Gamma^1_{21}$, 
$\Gamma^2_{21}$, $\Gamma^1_{22}$, $\Gamma^2_{22}$ from \mythetag{8.3}, due to 
\mythetag{8.4}, \mythetag{8.5}, \mythetag{8.6}, \mythetag{8.7}, \mythetag{8.8}, 
and \mythetag{8.9} they are expressed through scalar invariants 
$I_1$, $I_2$, $I_3$, $I_4$, $I_7$ and through the fields $N$ and $M$. According
to Theorem~\mythetheorem{8.2}, when expressing $j^{\,\sssize\text{Bgd}}_5$
we might need to differentiate these fields, i\.\,e\. calculate their covariant
derivatives along $\boldsymbol\alpha$ and $\boldsymbol\gamma$. Doing it, we shall
produce higher order invariants in the sequence given by \mythetag{8.11} and
some covariant derivatives of $N$ and $M$ along $\boldsymbol\alpha$ and 
$\boldsymbol\gamma$. Covariant derivatives of $N$ and $M$ are described
by the following lemma.
\mylemma{8.2} In the case where $\Omega=0$ covariant derivatives of the 
pseudoscalar fields $M$ and $N$ along the pseudovectorial fields 
$\boldsymbol\alpha$ and $\boldsymbol\gamma$ are expressed through the 
scalar invariants $I_1$, $I_2$, $I_3$, $I_4$ etc in the recurrent 
sequence given by \mythetag{8.11} and through these two fields themselves. 
\endproclaim
     Lemma~\mythelemma{8.2} in the present paper is similar to Lemma~8.2 in 
\mycite{12}. This lemma is proved by means of the following explicit formulas:
$$
\xalignat 2
&\hskip -2em
\nabla_{\boldsymbol\alpha}N=M,
&&\nabla_{\boldsymbol\gamma}N=0,\\
\vspace{-1.5ex}
\mytag{8.28}\\
\vspace{-1.5ex}
&\hskip -2em
\nabla_{\boldsymbol\alpha}M=I_4\,N^3+2\,I_1\,N\,M,
&&\nabla_{\boldsymbol\gamma}M=\sqrt{N^3\,I_7}\,N^2. 
\endxalignat
$$
The formulas \mythetag{8.28} are derived from the formulas \thetag{8.15} and
\thetag{8.17} in \mycite{12} by setting $\Omega=0$ in them. They should be 
applied repeatedly in order to calculate higher order covariant derivatives 
of $M$ and $N$.\par
     Now, combining Lemma~\mythelemma{8.2} with Theorem~\mythetheorem{8.2}
and with the arguments given just after Theorem~\mythetheorem{8.2}, we derive
the following theorem.
\mytheorem{8.3} Within the intersection class\/ {\rm ShrID1\,$\cap$\,BgdET4}
coinciding with {\rm BgdET4}, i\.\,e\. if the conditions 
\mythetag{7.3} are fulfilled, Bagderina's basic invariant 
$I^{\,\sssize\text{Bgd}}_2$ from \mythetag{6.4} can be expressed
through $I_1$, $I_3$ and through higher order invariants $I_4$,  
$I_6$, $I_7$ etc in the recurrent sequence given by \mythetag{8.11}. 
\endproclaim     
\head
9. Conclusions.
\endhead
     Three classes of umbilical equations, zero Gaussian curvature equations,
and zero mean curvature equations are defined in the present paper. They
specify the equations of the form \mythetag{1.1} in all cases of intermediate 
degeneration. Generally speaking, these geometric classes do not fit into 
particular subcases of both classifications from \myciterange{3}{3}{, }{4} 
and/or from \mycite{11}.\par
    Being intersected with the class {\rm ShrID1}, which corresponds to the 
first case of intermediate degeneration, the class of zero mean curvature 
equations produces a subclass coinciding with the intersection 
class\/ {\rm ShrID1\,$\cap$\,BgdET4} and with Bagderina's class {\rm BgdET4}
of type four equations. We have compared two classifications from 
\myciterange{3}{3}{, }{4} and from \mycite{11} within this intersection class. 
As a result we have found that most basic structures and basic formulas from 
Bagderina's paper \mycite{11} do coincide or are very closely related to those 
in \myciterange{3}{3}{, }{4}, though they are given in different notations (see Lemma~\mythelemma{7.1}, Lemma~\mythelemma{7.2}, and Lemma~\mythelemma{7.3}).
Similar results for the case of general position and for the other intersection
class\/ {\rm ShrID1\,$\cap$\,BgdET2} were obtained in \mycite{19} and
\mycite{12}.\par 
     For her type four equations class {\rm BgdET4} in \mycite{11} 
Yu\.\,Yu\.\,Bagderina introduces two basic invariants 
$I^{\,\sssize\text{Bgd}}_1$ and $I^{\,\sssize\text{Bgd}}_2$. For our class 
{\rm ShrID1}, which covers Bagderina's class {\rm BgdET4}, three basic invariants 
$I_1$, $I_2$, $I_3$ were introduced in \myciterange{3}{3}{, }{4}. However within 
the intersection class {\rm ShrID1\,$\cap$\,BgdET4} the invariant $I_2$ vanishes,
so we have two Bagderina's invariants versus two ours. In \mycite{20}
Yu\.\,Yu\.\,Bagderina managed to express the invariants $I_1$ and $I_3$
through her invariants $I^{\,\sssize\text{Bgd}}_1$, $I^{\,\sssize\text{Bgd}}_2$
and through their derivatives. In the present paper we present the converse 
result, i\.\,e\. we have proved that the invariants $I^{\,\sssize\text{Bgd}}_1$ 
and  $I^{\,\sssize\text{Bgd}}_2$ can be expressed through our invariants $I_1$, 
$I_2$ and through their derivatives (see Lemma~\mythelemma{8.1} and
Theorem~\mythetheorem{8.3}). Thus both sets of basic invariants are equally
applicable to solving the equivalence problem for the equations \mythetag{1.1}
within the intersection class {\rm ShrID1\,$\cap$\,BgdET4}.\par

\enddocument
\end